\documentclass[5p,times]{elsarticle}
\journal{Automatica}

\usepackage{amsfonts}
\usepackage{mathtools}
\usepackage{float}
\usepackage{amsmath}
\usepackage{graphicx}
\usepackage{amsthm}
\usepackage{epstopdf}
\usepackage{subfigure}
\usepackage{color}
\usepackage{epsfig}
\usepackage{multirow}
\usepackage{fancyhdr}
\usepackage{psfrag}
\usepackage{dsfont}
\usepackage{amssymb}
\usepackage{xcolor}
\usepackage[normalem]{ulem}
\usepackage{graphicx}
\usepackage{booktabs}
\usepackage{microtype}
\usepackage{mathrsfs}
\usepackage{moreverb}
\usepackage{graphicx}
\usepackage{float}
\usepackage{hyperref}
\usepackage{microtype}
\usepackage{bm}
\usepackage{amssymb}

\newtheorem{theorem}{Theorem}
\newtheorem{assumption}{Assumption}
\newtheorem{lemma}[theorem]{Lemma}
\newtheorem{define}[theorem]{Definition}

\newtheorem{proposition}[theorem]{Proposition}
\newtheorem{corollary}{Corollary}

\usepackage[figuresright]{rotating}
\biboptions{authoryear}

\begin{document}

\begin{frontmatter}

\title{Relative-energy stabilization of an irreversible port-Hamiltonian tubular reactor}
\tnotetext[mytitlenote]{}

\author[First]{Weijun. Zhou}
%\author[Second]{ Boussad. Hamroun}

\address[First]{School of Information and Electrical Engineering, Hangzhou City University, 310027, P. R. China,(e-mail: zhouwj@hzcu.edu.cn)}
%\address[Second]{LAGEPP, UMR CNRS 5007, University of Lyon1,  Villeurbanne, France, (e-mails: \{boussad.Hamroun;francoise.couenne\}@univ-lyon1.fr).}

\begin{abstract}

This paper develops a relative-energy control approach for a boundary-controlled irreversible port-Hamiltonian model of a non-isothermal tubular reactor around a spatially nonuniform stationary profile. The starting point is a thermodynamic duality between the entropy-based availability used in earlier reactor control and the relative Hamiltonian induced by the irreversible port-Hamiltonian representation. Under physically interpretable stability conditions---positive heat capacity, positive compressibility and positive compositional susceptibility---the Hamiltonian density is shown to be strictly convex, while the corresponding entropy representation is strictly concave. The resulting relative Hamiltonian decomposes into kinetic, compressive, thermal and mixing contributions and yields an exact relative-energy balance in which irreversible dissipation is separated from profile-dependent terms. This structure leads to a bounded monotone distributed jacket feedback that injects thermal damping without cancelling the nonlinear PDE dynamics. A continuous coercivity condition gives a local a-priori decay estimate for sufficiently regular nonlinear solutions, while a composite energy functional provides an exponential estimate for the linearized closed-loop initial-boundary value problem. Structure-preserving PFEM simulations are used to verify the predicted dissipation mechanisms and constrained profile regulation.

\end{abstract}


\end{frontmatter}

%\linenumbers

%===========================================================
\section{Introduction}
\label{sec:introduction}
%===========================================================

Tubular chemical reactors are distributed-parameter systems in which
convection, diffusion, heat transfer, and chemical reaction interact
over a spatial domain. Their control has been studied using
infinite-dimensional feedback, passivity-based methods, and model
predictive control
\citep{orlovDochain2002,franco2020,zhang2020}.
A complementary line of work exploits irreversible thermodynamics,
where entropy production and thermodynamic availability provide a
physically meaningful basis for stability and control
\citep{ydstieAlonso1997,alonsoYdstie2001,ydstie2002,ruszkowski2005}.
Availability-based control has been developed for stirred-tank
reactors \citep{hoang2012} and, more recently, for tubular reactors
\citep{zhou2021}.

The previous tubular-reactor formulation in \citet{zhou2021} assumes
constant density and pressure and prescribes the convection velocity.
Its availability function is therefore purely thermochemical and does
not contain momentum or compressive-energy contributions. The present work instead considers the BC-IPHS tubular-reactor model
developed in \citet{zhou2026bciphs}, which includes dynamic momentum,
variable density, physical inlet--outlet ports, reaction, diffusion,
heat conduction, and distributed jacket exchange. Its formulation is
consistent with the irreversible and boundary port-Hamiltonian framework
\citep{ramirez2013,ramirez2016,ramirez2022bciphs}. Port-Hamiltonian formulations are particularly useful in this setting
because they make the energy and boundary power exchanges explicit, while the irreversible extension also retains the entropy-production
structure \citep{vanderschaftMaschke2002,rashad2020,ramirezLeGorrec2022}. These two viewpoints are thermodynamically compatible rather than
independent. Availability exploits the concavity of entropy to measure
departure from a reference state, whereas the irreversible
port-Hamiltonian formulation organizes the corresponding energy
balance and entropy-production mechanisms. In the present reactor,
this common thermodynamic structure allows the availability to be
identified directly with the relative Hamiltonian of the BC-IPHS.

The main difficulty is that the desired reactor state is generally a
spatially nonuniform stationary profile rather than a thermodynamic
equilibrium. Reaction and transport fluxes may remain nonzero along
this profile, so the relative Hamiltonian does not reduce to a simple
incremental dissipation identity. Spatial variation of the reference
generates additional terms whose sign is not fixed a priori. The
control problem is therefore to identify a suitable relative storage,
separate the intrinsically dissipative mechanisms from the
profile-dependent terms, and determine whether bounded thermal damping
through the jacket is sufficient to obtain a local decay estimate.

The first contribution of this paper is to identify the
entropy-based thermodynamic availability with the entropy
representation of the relative Hamiltonian of the momentum-inclusive
BC-IPHS. For the adopted constitutive closure, strict positivity is
linked to explicit thermodynamic stability conditions, and the
resulting thermo-kinetic availability decomposes into kinetic,
compressive, thermal, and mixing contributions. This construction
extends the availability used in \citet{zhou2021} without claiming a
new general Bregman or Legendre duality, which is classical in
thermodynamics and relative-energy analysis
\citep{ramirez2016,feireisl2012}.

The second contribution is an exact relative-energy balance about the
spatially varying stationary profile. The mechanical subsystem is
shown to evolve autonomously under the adopted barotropic closure,
while the thermochemical balance explicitly separates reaction,
diffusion, heat conduction, reference-profile terms,
mechanical--thermochemical coupling, and jacket actuation. This
decomposition leads to a bounded monotone jacket feedback whose shifted thermal supply is nonpositive. A bulk--boundary coercivity
condition is then used to derive a conditional nonlinear
relative-energy estimate, while a pointwise Schur-complement condition
provides a convenient sufficient criterion for positive definiteness of the bulk matrix. For the linearized closed-loop system, a composite Lyapunov functional
combining thermochemical relative energy and a weighted mechanical
energy yields an exponential energy estimate for sufficiently regular
solutions of the linearized initial-boundary-value problem.

The analysis is evaluated using the structure-preserving PFEM model of
the reactor. The numerical results confirm that the selected gain satisfies the
derived pointwise bulk-positivity condition and examine the associated
semi-discrete coercivity, together with the bounded jacket action and
the acceleration of the distributed temperature--composition transient. A temporary throughput
perturbation is also considered to activate the mechanical dynamics
and illustrate the one-way coupling from the mechanical subsystem
to the thermochemical dynamics. The selected
open-loop operating point is already numerically convergent; the role
of the controller is therefore to increase structured thermal
dissipation and improve the transient rather than to stabilize an
artificially unstable benchmark.

The remainder of the paper is organized as follows.
Section~\ref{sec:model} recalls the BC-IPHS reactor model.
Section~\ref{sec:TKA} introduces the thermo-kinetic availability and
its relative-Hamiltonian interpretation.
Section~\ref{sec:relative_balance} derives the exact relative-energy
balance.
Section~\ref{sec:control} presents the bounded jacket feedback and the
stability estimates.
Section~\ref{sec:numerical} gives the numerical verification, followed
by the conclusions.

%===========================================================
\section{Boundary-Controlled Irreversible Port-Hamiltonian Reactor}
\label{sec:model}
%===========================================================

The control design is based on the boundary-controlled irreversible
port-Hamiltonian (BC-IPHS) tubular-reactor model developed in
\citep{zhou2026bciphs}. The complete port-Hamiltonian realization,
boundary ports, and structure-preserving PFEM discretization are given
there. Only the model properties required for the present stability
analysis are recalled below.

%-----------------------------------------------------------
\subsection{Reactor model and thermodynamic structure}
\label{subsec:model_structure}
%-----------------------------------------------------------

Consider a one-dimensional binary reacting mixture
$A\rightleftharpoons B$ on $z\in[0,L]$. The conservative state and
derived variables are
\[
 x=
 \begin{bmatrix}
 \rho_A & \rho_B & m & \eta
 \end{bmatrix}^{\top},
 \qquad
 \rho=\rho_A+\rho_B,
 \qquad
 \alpha=\frac{\rho_A}{\rho},
 \qquad
 v=\frac{m}{\rho},
 \qquad
 s=\frac{\eta}{\rho},
\]
on the admissible domain
$\rho>0$, $T>0$, and $0<\alpha<1$.

The Hamiltonian density is
\begin{equation}
 h(x)
 =
 \frac{m^2}{2\rho}
 +
 \Pi(\rho)
 +
 \rho u_{tc}(s,\alpha),
 \label{eq:hamiltonian-density}
\end{equation}
with total Hamiltonian
$H[x]=\int_0^L h(x)\,dz$.
Here $\Pi(\rho)=\rho a_c(\rho)$ is the compressive energy density and
$u_{tc}$ is the thermochemical specific internal energy.

Let
\[
 \varepsilon(\rho_A,\rho_B,\eta)
 =
 \Pi(\rho)+\rho u_{tc}(s,\alpha).
\]
The density-form Gibbs relation is
\[
 d\varepsilon
 =
 T\,d\eta
 +
 \mu_A\,d\rho_A
 +
 \mu_B\,d\rho_B,
\]
while separation of the compressive contribution gives
\begin{equation}
 du_{tc}
 =
 T\,ds+A_r\,d\alpha,
 \qquad
 A_r:=\mu_A-\mu_B.
 \label{eq:gibbs}
\end{equation}
The common compressive contribution to the chemical potentials
therefore cancels from the affinity. The Hamiltonian efforts are
\begin{equation}
 \nabla_x h
 =
 \begin{bmatrix}
 \mu_A-v^2/2\\
 \mu_B-v^2/2\\
 v\\
 T
 \end{bmatrix}.
 \label{eq:model_hamiltonian_gradient}
\end{equation}

The pressure is generated by the compressive energy according to
\begin{equation}
 p(\rho)
 =
 \rho\Pi'(\rho)-\Pi(\rho).
 \label{eq:model_pressure}
\end{equation}
For the weakly compressible benchmark,
\[
 p(\rho)
 =
 p_{\rm ref}
 +
 K_{\rm bulk}
 \left(
 \frac{\rho}{\rho_{\rm ref}}-1
 \right),
 \qquad
 p'(\rho)
 =
 \frac{K_{\rm bulk}}{\rho_{\rm ref}}>0.
\]
The caloric ideal-mixture closure further gives
\[
 c_p(\alpha)
 =
 \alpha c_{p,A}
 +(1-\alpha)c_{p,B}>0,
 \qquad
 A_\alpha
 :=
 \left.
 \frac{\partial A_r}{\partial\alpha}
 \right|_T
 =
 \frac{R_sT}{\alpha(1-\alpha)}>0.
\]
The conditions $p'(\rho)>0$, $c_p(\alpha)>0$, and $A_\alpha>0$ are
used below to establish strict convexity of the relative storage.

The species balances are
\begin{align*}
 \partial_t\rho_A+\partial_z(\rho_Av+J)&=-r,\\
 \partial_t\rho_B+\partial_z(\rho_Bv-J)&=r.
\end{align*}
Their sum and the momentum equation give the autonomous mechanical
subsystem
\begin{align}
 \partial_t\rho+\partial_zm
 &=0,
 \label{eq:model_mass}\\
 \partial_tm
 +
 \partial_z
 \left(
 \frac{m^2}{\rho}+p(\rho)
 \right)
 &=0.
 \label{eq:model_momentum}
\end{align}

The entropy-conjugate thermodynamic forces are
\begin{equation}
 X_r=\frac{A_r}{T},
 \qquad
 X_d=-\frac{1}{T}\partial_zA_r,
 \qquad
 X_q=\partial_z\left(\frac{1}{T}\right).
 \label{eq:model_forces}
\end{equation}
For the constitutive closure considered here,
\begin{equation}
 r
 =
 k_f(T)\rho_A
 \left[
 1-\exp\left(-\frac{A_r}{R_sT}\right)
 \right],
 \qquad
 k_f(T)=k_0\exp\left(-\frac{E_a}{RT}\right),
 \label{eq:model_reaction}
\end{equation}
and
\begin{equation}
 J=-M_d\partial_zA_r,
 \qquad
 M_d=\frac{D\rho\alpha(1-\alpha)}{R_sT},
 \qquad
 q=-k\partial_zT.
 \label{eq:model_transport}
\end{equation}
These laws satisfy
\[
 rX_r\ge0,
 \qquad
 JX_d\ge0,
 \qquad
 qX_q\ge0.
\]

The jacket interaction is
\begin{equation}
 q_j=\kappa_j(T_j-T),
 \qquad
 f_j:=\frac{q_j}{T},
 \label{eq:model_jacket}
\end{equation}
where $T_j(z,t)>0$ is the distributed control input. The entropy
balance is
\begin{equation}
 \partial_t(\rho s)
 +
 \partial_z
 \left(
 \rho sv+\frac{q}{T}
 \right)
 =
 \sigma_r+\sigma_d+\sigma_q+\frac{q_j}{T},
 \label{eq:entropy-balance}
\end{equation}
where
\[
 \sigma_r=\frac{A_rr}{T},
 \qquad
 \sigma_d=-\frac{J}{T}\partial_zA_r,
 \qquad
 \sigma_q=q\partial_z\left(\frac1T\right)
\]
are nonnegative. The jacket exchange similarly satisfies
\[
 \frac{q_j}{T}
 =
 \frac{q_j}{T_j}+\sigma_j,
 \qquad
 \sigma_j
 =
 \kappa_j\frac{(T_j-T)^2}{TT_j}
 \ge0.
\]

With
\[
 h_{\rm th}
 =
 a_c(\rho)+u_{tc}(s,\alpha)+\frac{p(\rho)}{\rho},
\]
the axial total-energy flux is
\begin{equation}
 F_H
 =
 \rho v
 \left(
 h_{\rm th}+\frac{v^2}{2}
 \right)
 +
 A_rJ+q.
 \label{eq:model_energy_flux}
\end{equation}
The corresponding first-law balance is
\[
 \frac{dH}{dt}
 =
 -[F_H]_0^L
 +
 \int_0^L q_j\,dz.
\]

%-----------------------------------------------------------
\subsection{Boundary conditions and stationary profile}
\label{subsec:model_boundary}
%-----------------------------------------------------------

The analysis is restricted to positive subsonic flow,
\begin{equation}
 0<v<c(\rho),
 \qquad
 c^2(\rho):=p'(\rho).
 \label{eq:model_subsonic}
\end{equation}
At the inlet, the total species and energy fluxes are prescribed,
\[
 (\rho_Av+J)(0,t)=F_A^{\rm in},
 \qquad
 (\rho_Bv-J)(0,t)=F_B^{\rm in},
 \qquad
 F_H(0,t)=F_H^{\rm in},
\]
whereas the outlet conditions are
\begin{equation}
 J(L,t)=0,
 \qquad
 q(L,t)=0,
 \qquad
 p(\rho(L,t))=p_{\rm out}.
 \label{eq:model_bc_out}
\end{equation}
The inlet species conditions imply
\begin{equation}
 m(0,t)
 =
 F_A^{\rm in}+F_B^{\rm in}
 =:m_d.
 \label{eq:model_bc_mass}
\end{equation}

Let
\[
 x_d(z)
 =
 \begin{bmatrix}
 \rho_{A,d}(z)&
 \rho_{B,d}(z)&
 m_d&
 \eta_d(z)
 \end{bmatrix}^{\top}
\]
be a sufficiently smooth stationary solution corresponding to the
stationary jacket profile $T_{j,d}(z)$ and the same axial boundary
data.

The stationary mass and momentum balances imply
\[
 \partial_zm_d=0,
 \qquad
 \partial_z
 \left(
 \frac{m_d^2}{\rho_d}+p(\rho_d)
 \right)=0.
\]
Moreover,
\[
 \frac{d}{d\rho}
 \left(
 \frac{m_d^2}{\rho}+p(\rho)
 \right)
 =
 c^2(\rho)-v^2>0
\]
on the positive subsonic branch. Hence the stationary mechanical
variables are spatially constant:
\begin{equation}
 \rho_d=\rho_\star,
 \qquad
 m_d=m_\star,
 \qquad
 v_d=v_\star.
 \label{eq:constant-mechanics}
\end{equation}
In contrast, $T_d(z)$ and $\alpha_d(z)$ are generally nonuniform, and
the stationary reaction and transport quantities $r_d$, $J_d$, and
$q_d$ need not vanish. This distinction is central to the
relative-energy balance derived below.

%===========================================================
\section{Thermo-Kinetic Availability and Relative Hamiltonian}
\label{sec:TKA}
%===========================================================

This section constructs the relative storage used in the control
analysis. Both constructions originate from the same thermodynamic
energy--entropy structure, but emphasize different representations
of it. The entropy--energy tangent-gap duality itself is classical;
the purpose here is to identify the thermodynamic availability of the
reactor with the relative Hamiltonian of the momentum-inclusive
BC-IPHS, establish its positivity, and determine its physical
decomposition for the constitutive closure considered below.

%-----------------------------------------------------------
\subsection{Entropy representation and thermo-kinetic availability}
\label{subsec:TKA_definition}
%-----------------------------------------------------------

Let $\mathcal E=h(x)$ and introduce the entropy representation
\begin{equation}
 \eta=\mathscr S(\mathcal E,\rho_A,\rho_B,m),
 \qquad
 y=
 \begin{bmatrix}
 \mathcal E & \rho_A & \rho_B & m
 \end{bmatrix}^{\top}.
 \label{eq:entropy_representation}
\end{equation}
Define the species efforts
\[
 e_A:=\mu_A-\frac{v^2}{2},
 \qquad
 e_B:=\mu_B-\frac{v^2}{2}.
\]
Using \eqref{eq:model_hamiltonian_gradient},
\[
 d\mathcal E
 =
 e_A\,d\rho_A
 +
 e_B\,d\rho_B
 +
 v\,dm
 +
 T\,d\eta,
\]
and therefore
\begin{equation}
 \boxed{
 \nabla_y\mathscr S(y)
 =
 \begin{bmatrix}
 1/T &
 -e_A/T &
 -e_B/T &
 -v/T
 \end{bmatrix}^{\top}.
 }
 \label{eq:entropy_gradient}
\end{equation}
Thus the momentum contribution $-v/T$ follows directly from the
energy--entropy transformation rather than being introduced as an
additional mechanical correction.

Let $x_d(z)$ denote the stationary reference profile and let $y_d(z)$
be its entropy representation.

\begin{define}[Thermo-kinetic availability]
\label{def:TKA}
The entropy-form thermo-kinetic availability density relative to
$y_d$ is
\begin{equation}
 \boxed{
 \begin{aligned}
 \bar a_{\rm TK}(y\mid y_d)
 &=
 \mathscr S(y_d)
 +
 \nabla_y\mathscr S(y_d)^{\top}(y-y_d)
 -
 \mathscr S(y),
 \\
 a_{\rm TK}(y\mid y_d)
 &=
 T_d\,\bar a_{\rm TK}(y\mid y_d),
 \qquad T_d>0 .
 \end{aligned}
 }
 \label{eq:TKA_definition}
\end{equation}
\end{define}

Here ``kinetic'' refers to the macroscopic kinetic energy associated
with the momentum state $m$, not to chemical-reaction kinetics.

%-----------------------------------------------------------
\subsection{Duality, convexity, and positivity}
\label{subsec:TKA_convexity}
%-----------------------------------------------------------

\begin{proposition}[Entropy--Hamiltonian duality]
\label{prop:TKA_duality}
The thermo-kinetic availability satisfies
\begin{equation}
 \boxed{
 a_{\rm TK}(y\mid y_d)
 =
 D_h(x\mid x_d)
 =
 T_dD_{-\mathscr S}(y\mid y_d),
 }
 \label{eq:TKA_duality}
\end{equation}
where
\begin{equation}
 D_h(x\mid x_d)
 :=
 h(x)-h(x_d)
 -
 \nabla h(x_d)^{\top}(x-x_d)
 \label{eq:relative_Hamiltonian}
\end{equation}
and
\[
 D_{-\mathscr S}(y\mid y_d)
 :=
 -\mathscr S(y)+\mathscr S(y_d)
 +
 \nabla\mathscr S(y_d)^{\top}(y-y_d).
\]
\end{proposition}

\begin{proof}
Using \eqref{eq:entropy_gradient} at the reference state and
$\mathscr S(y)=\eta$ gives
\begin{align*}
 T_dD_{-\mathscr S}(y\mid y_d)
 ={}&
 \mathcal E-\mathcal E_d
 -e_{A,d}(\rho_A-\rho_{A,d})
 \\
 &-
 e_{B,d}(\rho_B-\rho_{B,d})
 -v_d(m-m_d)
 -T_d(\eta-\eta_d).
\end{align*}
Since $\mathcal E=h(x)$ and
\[
 \nabla h(x_d)
 =
 \begin{bmatrix}
 e_{A,d} &
 e_{B,d} &
 v_d &
 T_d
 \end{bmatrix}^{\top},
\]
the right-hand side is exactly $D_h(x\mid x_d)$.
Definition~\ref{def:TKA} then gives
\eqref{eq:TKA_duality}.
\end{proof}

\begin{proposition}[Strict convexity of the Hamiltonian]
\label{prop:Hamiltonian_convexity}
Assume, on the admissible state domain,
\[
 \rho>0,\qquad
 T>0,\qquad
 c_p>0,\qquad
 p'(\rho)>0,\qquad
 A_\alpha
 :=
 \left.
 \frac{\partial A_r}{\partial\alpha}
 \right|_T
 >0.
\]
Then $u_{tc}$ is strictly convex in $(s,\alpha)$ and the Hamiltonian
density \eqref{eq:hamiltonian-density} is strictly convex.
\end{proposition}

\begin{proof}
Let
\[
 \zeta=
 \begin{bmatrix}
 s\\
 \alpha
 \end{bmatrix},
 \qquad
 A_T:=
 \left.
 \frac{\partial A_r}{\partial T}
 \right|_\alpha .
\]
Since
$T=\partial u_{tc}/\partial s$,
$A_r=\partial u_{tc}/\partial\alpha$, and
$T_s=T/c_p$, symmetry of the mixed derivatives gives
\[
 H_{tc}
 :=
 \nabla_\zeta^2u_{tc}
 =
 \begin{bmatrix}
 \displaystyle\frac{T}{c_p}
 &
 \displaystyle\frac{T}{c_p}A_T
 \\[2mm]
 \displaystyle\frac{T}{c_p}A_T
 &
 \displaystyle
 A_\alpha+\frac{T}{c_p}A_T^2
 \end{bmatrix}.
\]
Its leading principal minor and determinant satisfy
\[
 \frac{T}{c_p}>0,
 \qquad
 \det H_{tc}
 =
 \frac{T}{c_p}A_\alpha>0,
\]
hence $H_{tc}\succ0$.

For an arbitrary perturbation of the complete state,
with $\delta\rho=\delta\rho_A+\delta\rho_B$,
\begin{equation}
 \boxed{
 \delta^2h
 =
 \frac{1}{\rho}
 \left(
 \delta m-v\,\delta\rho
 \right)^2
 +
 \frac{p'(\rho)}{\rho}
 (\delta\rho)^2
 +
 \rho\,
 \delta\zeta^{\top}
 H_{tc}
 \delta\zeta .
 }
 \label{eq:full_hessian_quadratic}
\end{equation}
If the right-hand side vanishes, $p'(\rho)>0$ gives
$\delta\rho=0$, the kinetic term then gives $\delta m=0$, and
$H_{tc}\succ0$ gives $\delta\zeta=0$. Hence $h$ is strictly convex.
\end{proof}

The corresponding Hessian relation in entropy variables follows by differentiating twice
the implicit relation
\[
 h\!\left(
 \rho_A,\rho_B,m,\mathscr S(y)
 \right)
 =
 \mathcal E .
\]
For corresponding perturbations,
\[
 \delta y^{\top}
 \nabla_y^2\mathscr S(y)
 \delta y
 =
 -\frac{1}{T}
 \delta x^{\top}
 \nabla_x^2h(x)
 \delta x .
\]
Since $T>0$, the energy--entropy transformation is locally invertible,
and Proposition~\ref{prop:Hamiltonian_convexity} therefore implies
strict concavity of $\mathscr S$ in the present fixed-volume density
representation. This does not contradict the usual homogeneity
degeneracy of a fully extensive entropy, because the present variables
are densities and the compressive contribution with $p'(\rho)>0$
removes the free radial scaling direction.

Combining Propositions~\ref{prop:TKA_duality} and
\ref{prop:Hamiltonian_convexity} gives, in any convex admissible
neighborhood of the reference profile,
\begin{equation}
 \boxed{
 a_{\rm TK}(y\mid y_d)\ge0,
 \qquad
 a_{\rm TK}(y\mid y_d)=0
 \quad\Longleftrightarrow\quad
 x=x_d .
 }
 \label{eq:TKA_positive}
\end{equation}

\begin{figure}[t]
 \centering
 \includegraphics[width=\columnwidth]
 {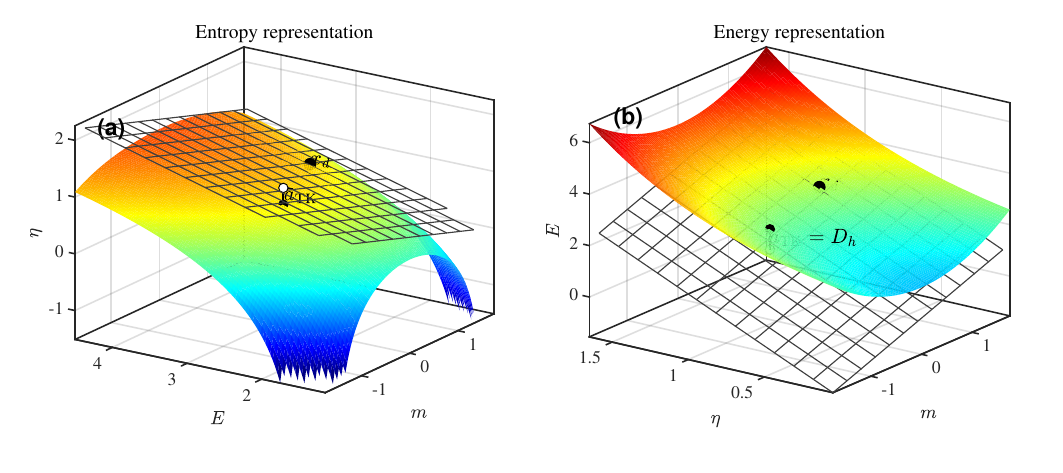}
 \caption{Dual geometry of the thermo-kinetic availability.
 (a) Entropy representation: $\bar a_{\rm TK}$ is the tangent gap
 associated with the concave entropy surface.
 (b) Energy representation:
 $a_{\rm TK}=D_h=T_d\bar a_{\rm TK}$ is the corresponding gap between
 the convex energy surface and its supporting tangent plane.}
 \label{fig:TKA_geometry}
\end{figure}

The dual geometry is illustrated in
Fig.~\ref{fig:TKA_geometry}. In the entropy representation,
$\bar a_{\rm TK}$ is the gap between the supporting tangent plane and
the concave entropy surface. In the energy representation, the
corresponding energy gap is
$a_{\rm TK}=D_h=T_d\bar a_{\rm TK}$. The figure is only a normalized
illustration of the duality and does not represent the reactor
constitutive surface.

%-----------------------------------------------------------
\subsection{Exact decomposition and relation to the previous availability}
\label{subsec:TKA_decomposition}
%-----------------------------------------------------------

For the constitutive closure adopted in this paper, the relative
Hamiltonian admits a physical decomposition.

\begin{proposition}[Decomposition of the relative Hamiltonian]
\label{prop:TKA_decomposition}
For the caloric ideal-mixture closure and the barotropic pressure law
of Section~\ref{sec:model},
\begin{equation}
 \boxed{
 a_{\rm TK}
 =
 a_{\rm kin}
 +
 a_{\rm comp}
 +
 a_{\rm th}
 +
 a_{\rm mix}.
 }
 \label{eq:TKA_decomposition}
\end{equation}
The four contributions are
\begin{align*}
 a_{\rm kin}
 &=
 \frac12\rho(v-v_d)^2,
 \\[1mm]
 a_{\rm comp}
 &=
 \frac{K_{\rm bulk}}{\rho_{\rm ref}}
 \left[
 \rho\ln\frac{\rho}{\rho_d}
 -\rho+\rho_d
 \right],
 \\[1mm]
 a_{\rm th}
 &=
 \rho c_p(\alpha)
 \left[
 T-T_d-T_d\ln\frac{T}{T_d}
 \right],
 \\[1mm]
 a_{\rm mix}
 &=
 \rho R_sT_d
 \left[
 \alpha\ln\frac{\alpha}{\alpha_d}
 +(1-\alpha)
 \ln\frac{1-\alpha}{1-\alpha_d}
 \right].
\end{align*}
\end{proposition}

\begin{proof}
The kinetic Bregman remainder satisfies
\[
 \frac{m^2}{2\rho}
 -\frac{m_d^2}{2\rho_d}
 +\frac{v_d^2}{2}(\rho-\rho_d)
 -v_d(m-m_d)
 =
 \frac12\rho(v-v_d)^2.
\]

For the compressive contribution,
\[
 D_\Pi(\rho\mid\rho_d)
 =
 \Pi(\rho)-\Pi(\rho_d)
 -\Pi'(\rho_d)(\rho-\rho_d).
\]
Since
\[
 \Pi''(\rho)
 =
 \frac{p'(\rho)}{\rho}
 =
 \frac{K_{\rm bulk}}
 {\rho_{\rm ref}\rho},
\]
direct integration gives the expression for $a_{\rm comp}$ above.

For the thermochemical contribution,
\[
 D_{tc}
 =
 \rho
 \left[
 u_{tc}(s,\alpha)
 -u_{tc}(s_d,\alpha_d)
 -T_d(s-s_d)
 -A_{r,d}(\alpha-\alpha_d)
 \right].
\]
Substitution of the caloric ideal-mixture energy and entropy relations
gives
\[
 \frac{D_{tc}}{\rho}
 =
 c_p(\alpha)
 \left[
 T-T_d-T_d\ln\frac{T}{T_d}
 \right]
 +
 R_sT_d
 \left[
 \alpha\ln\frac{\alpha}{\alpha_d}
 +(1-\alpha)
 \ln\frac{1-\alpha}{1-\alpha_d}
 \right].
\]
Combining the kinetic, compressive, and thermochemical remainders
proves \eqref{eq:TKA_decomposition}.
\end{proof}

The decomposition shows explicitly that the momentum-inclusive model
adds kinetic and compressive storage directions to the thermal and
compositional thermodynamic availability.

On the fixed-density, prescribed-velocity manifold of the model used
in \citet{zhou2021},
\[
 a_{\rm kin}=a_{\rm comp}=0,
\]
so that $a_{\rm TK}/T_d$ recovers the full entropy-based thermodynamic
availability of that formulation. Thus the present construction
extends, rather than replaces, the previous availability to the
momentum-inclusive BC-IPHS state space.

The reduced availability introduced in \citet{zhou2021} retains only
the thermal contribution. With the notation of that work,
\begin{equation}
 \boxed{
 a_{\rm th}=T_d a_r.
 }
 \label{eq:JPC_recovery}
\end{equation}
The compositional term $a_{\rm mix}$ is absent from that reduced
storage by construction.

%===========================================================
\section{Relative-Energy Balance Along the Reactor Dynamics}
\label{sec:relative_balance}
%===========================================================

Section~\ref{sec:TKA} identified the thermo-kinetic availability with
the relative Hamiltonian of the BC-IPHS. We now derive its exact
evolution along the reactor dynamics.

Let $x_d(z)$ be a sufficiently smooth stationary solution with the
same axial boundary data. As established in
\eqref{eq:constant-mechanics}, its mechanical variables
$\rho_d$, $m_d$, and $v_d$ are spatially constant, whereas
$T_d(z)$ and $\alpha_d(z)$ may vary along the reactor. The stationary
reaction and transport fluxes therefore need not vanish. For any
quantity $\chi$, write $\delta\chi:=\chi-\chi_d$.

%-----------------------------------------------------------
\subsection{Mechanical relative-energy balance}
\label{subsec:mechanical_relative_balance}
%-----------------------------------------------------------

Under the adopted barotropic closure, the mass--momentum subsystem is
autonomous:
\[
 \partial_t\rho+\partial_zm=0,
 \qquad
 \partial_tm+
 \partial_z\left(\frac{m^2}{\rho}+p(\rho)\right)=0.
\]
Its mechanical energy density and flux are
\[
 h_m(\rho,m)
 =
 \frac{m^2}{2\rho}+\Pi(\rho),
 \qquad
 F_m
 =
 v\left[
 \frac12\rho v^2+\Pi(\rho)+p(\rho)
 \right],
\]
with $p(\rho)=\rho\Pi'(\rho)-\Pi(\rho)$. They satisfy
$\partial_t h_m+\partial_zF_m=0$.

The corresponding mechanical relative-energy density is
\begin{equation}
 a_m
 =
 \frac12\rho(v-v_d)^2
 +
 \Pi(\rho)-\Pi(\rho_d)
 -
 \Pi'(\rho_d)(\rho-\rho_d).
 \label{eq:mechanical_relative_explicit}
\end{equation}

\begin{proposition}[Mechanical relative-energy balance and boundary dissipation]
\label{prop:mechanical_boundary}
The density \eqref{eq:mechanical_relative_explicit} satisfies
\begin{equation}
 \partial_ta_m+\partial_zQ_m^{\rm rel}=0,
 \label{eq:mechanical_relative_identity}
\end{equation}
where
\begin{align}
 Q_m^{\rm rel}
 ={}&
 F_m-F_{m,d}
 -
 \left(
 \Pi'(\rho_d)-\frac12v_d^2
 \right)(m-m_d)
 \nonumber\\
 &-
 v_d
 \left[
 \frac{m^2}{\rho}+p(\rho)
 -
 \left(
 \frac{m_d^2}{\rho_d}+p_d
 \right)
 \right].
 \label{eq:mechanical_relative_flux}
\end{align}

If the solution remains on the positive subsonic branch and satisfies
the physical mechanical boundary conditions
\[
 m(0,t)=m_d,
 \qquad
 p(\rho(L,t))=p(\rho_d),
\]
then
\begin{equation}
 Q_m^{\rm rel}(0,t)\le0,
 \qquad
 Q_m^{\rm rel}(L,t)\ge0.
 \label{eq:mechanical_boundary_sign}
\end{equation}
Consequently,
\begin{equation}
 \boxed{
 \frac{d}{dt}\mathcal A_m
 =
 Q_m^{\rm rel}(0,t)-Q_m^{\rm rel}(L,t)
 \le0,
 \qquad
 \mathcal A_m:=\int_0^L a_m\,dz .
 }
 \label{eq:mechanical_storage_nonincrease}
\end{equation}
\end{proposition}

\begin{proof}
Since the reference mechanical state is spatially constant,
\[
 \nabla h_m(\rho_d,m_d)
 =
 \begin{bmatrix}
 \Pi'(\rho_d)-v_d^2/2\\
 v_d
 \end{bmatrix}
\]
is constant. Subtracting the corresponding linear combination of the
mass and momentum conservation laws from the mechanical energy balance
therefore gives
\eqref{eq:mechanical_relative_identity} with the relative flux
\eqref{eq:mechanical_relative_flux}.

At the inlet, $m=m_d$ and hence $v=m_d/\rho$. Along this boundary
manifold,
\[
 \frac{dQ_m^{\rm rel}}{d\rho}
 =
 m_d
 \bigl[c^2(\rho)-v^2\bigr]
 \left(
 \frac1\rho-\frac1{\rho_d}
 \right).
\]
Because $m_d>0$, $c^2(\rho)-v^2>0$, and
$Q_m^{\rm rel}(\rho_d)=0$, integration from $\rho_d$ to the current
inlet density yields
$Q_m^{\rm rel}(0,t)\le0$.

At the outlet, $p'(\rho)>0$ makes $p$ strictly increasing on the
admissible density interval. Hence
$p(\rho(L,t))=p(\rho_d)$ implies $\rho(L,t)=\rho_d$, and substitution
into \eqref{eq:mechanical_relative_flux} gives
\[
 Q_m^{\rm rel}(L,t)
 =
 \frac12\rho_d v(L,t)
 \bigl[v(L,t)-v_d\bigr]^2
 \ge0,
\]
where $v(L,t)>0$ follows from the positive-flow assumption.
Integration of \eqref{eq:mechanical_relative_identity} over
$[0,L]$ then gives \eqref{eq:mechanical_storage_nonincrease}.
\end{proof}

Thus, under the present barotropic closure, active jacket feedback is
not required to obtain non-increase of the mechanical relative energy.
The subsystem $(\rho,m)$ evolves autonomously, and its relative-energy
exchange occurs through the physical inlet and outlet.

%-----------------------------------------------------------
\subsection{Thermochemical relative-energy balance}
\label{subsec:thermochemical_relative_balance}
%-----------------------------------------------------------

Introduce
\[
 \zeta=
 \begin{bmatrix}
 s\\
 \alpha
 \end{bmatrix},
 \qquad
 w=\nabla_\zeta u_{tc}
 =
 \begin{bmatrix}
 T\\
 A_r
 \end{bmatrix},
 \qquad
 g:=\frac{q}{T},
\]
and define the specific thermochemical relative energy
\begin{equation}
 \psi(\zeta\mid\zeta_d)
 =
 u_{tc}(\zeta)
 -
 u_{tc}(\zeta_d)
 -
 w_d^\top(\zeta-\zeta_d),
 \qquad
 a_{tc}=\rho\psi .
 \label{eq:psi_definition}
\end{equation}

Using
$D_t:=\partial_t+v\partial_z$ and
\[
 \Sigma:=rX_r+JX_d+qX_q,
\]
the composition and entropy balances take the material form
\[
 \rho D_t\alpha=-r-\partial_zJ,
 \qquad
 \rho D_ts=\Sigma+f_j-\partial_zg.
\]
The corresponding stationary equations are
\[
 m_d\alpha_d'=-r_d-J_d',
 \qquad
 m_ds_d'=\Sigma_d+f_{j,d}-g_d'.
\]

Let
\[
 H_d:=\nabla_\zeta^2u_{tc}(\zeta_d),
 \qquad
 r_w:=\delta w-H_d\delta\zeta .
\]
If $u_{tc}\in C^3$ in a neighborhood of the reference profile, then
$r_w=O(|\delta\zeta|^2)$ locally; more precisely,
$|r_w|\le C|\delta\zeta|^2$ after restricting the neighborhood if
necessary.

\begin{proposition}[Thermochemical relative-energy identity]
\label{prop:thermochemical_identity}
For every sufficiently regular solution with $T>0$,
\begin{equation}
 \boxed{
 \partial_t(\rho\psi)
 +
 \partial_zQ_{tc}^{\rm rel}
 =
 \delta T\,\delta f_j
 +
 R_r+R_d+R_q+R_{\rm ref}+R_{\rm mt}.
 }
 \label{eq:exact_tc_identity}
\end{equation}
The relative thermochemical flux is
\begin{equation}
 Q_{tc}^{\rm rel}
 =
 m\psi
 +
 \delta A_r\,\delta J
 +
 \delta T\,\delta g,
 \label{eq:tc_relative_flux}
\end{equation}
and the remainder terms are
\begin{equation}
\begin{aligned}
 R_r
 &=
 -T_d\,\delta r\,\delta X_r
 +
 r_d\,\delta T\,\delta X_r,
 \\
 R_d
 &=
 -T_d\,\delta J\,\delta X_d
 +
 J_d\,\delta T\,\delta X_d,
 \\
 R_q
 &=
 \delta T
 \left(
 qX_q-q_dX_{q,d}
 \right)
 +
 \delta T_z\,\delta g,
 \\
 R_{\rm ref}
 &=
 m_d r_w^\top\zeta_d',
 \\
 R_{\rm mt}
 &=
 -\delta m\,
 (w_d')^\top\delta\zeta .
\end{aligned}
\label{eq:tc_remainders}
\end{equation}
\end{proposition}

\begin{proof}
Differentiating \eqref{eq:psi_definition} with respect to the current
state gives
$\nabla_\zeta\psi=\delta w$.
Since the reference profile depends on $z$,
\[
 \left.\partial_z\psi\right|_{\zeta}
 =
 -(w_d')^\top\delta\zeta,
\]
and hence
\[
 D_t\psi
 =
 \delta w^\top D_t\zeta
 -
 v(w_d')^\top\delta\zeta.
\]
Using the continuity equation therefore yields
\[
 \partial_t(\rho\psi)
 +
 \partial_z(m\psi)
 =
 \delta T(\Sigma+f_j-g_z)
 -
 \delta A_r(r+J_z)
 -
 m(w_d')^\top\delta\zeta.
\]

Adding
$\partial_z(\delta T\,\delta g)$ and
$\partial_z(\delta A_r\,\delta J)$
to the left-hand side and using the stationary equations gives
\begin{align*}
 \partial_t(\rho\psi)
 +
 \partial_zQ_{tc}^{\rm rel}
 ={}&
 \delta T\,\delta f_j
 +
 \delta T(\Sigma-\Sigma_d)
 -
 \delta A_r\,\delta r
 \\
 &+
 \delta T_z\,\delta g
 +
 \delta A_{r,z}\,\delta J
 \\
 &+
 m_d\delta w^\top\zeta_d'
 -
 m(w_d')^\top\delta\zeta .
\end{align*}

Because $w_d'=H_d\zeta_d'$, the last line can be written as
\[
 m_d
 \left(
 \delta w-H_d\delta\zeta
 \right)^\top\zeta_d'
 -
 \delta m\,
 (w_d')^\top\delta\zeta
 =
 R_{\rm ref}+R_{\rm mt}.
\]

It remains to separate the irreversible contributions. Using
$A_r=TX_r$ gives
\[
 \delta T
 \left(
 rX_r-r_dX_{r,d}
 \right)
 -
 \delta A_r\,\delta r
 =
 R_r,
\]
whereas
$A_{r,z}=-TX_d$ gives
\[
 \delta T
 \left(
 JX_d-J_dX_{d,d}
 \right)
 +
 \delta A_{r,z}\,\delta J
 =
 R_d.
\]
The remaining heat-conduction terms give $R_q$ in
\eqref{eq:tc_remainders}, completing the proof.
\end{proof}

For the constant Fourier conductivity $k>0$ adopted here, the
heat-conduction remainder has a useful exact form. Let
\[
 Y:=\frac{T}{T_d}.
\]
Since
\[
 \frac{q}{T}
 =
 -k
 \left(
 \frac{Y_z}{Y}
 +
 \frac{T_d'}{T_d}
 \right),
 \qquad
 \delta g=-k\frac{Y_z}{Y},
\]
substitution into the expression for $R_q$ gives
\begin{equation}
 \boxed{
 R_q
 =
 -k\frac{T_d}{Y^2}Y_z^2
 +
 kT_d'
 \left(
 1-\frac{1}{Y}
 \right)Y_z .
 }
 \label{eq:Rq_ratio}
\end{equation}

The first term in \eqref{eq:Rq_ratio} is nonpositive, whereas the
second is generated by the spatial variation of the stationary
temperature profile and has no fixed sign. Likewise, the terms
proportional to $r_d$ and $J_d$ in \eqref{eq:tc_remainders} arise
because the stationary reactor may sustain nonzero reaction and
diffusive fluxes. These profile-dependent terms are therefore retained
explicitly rather than assigned a sign at this stage.

%-----------------------------------------------------------
\subsection{Distributed thermo-kinetic availability balance}
\label{subsec:distributed_TKA_balance}
%-----------------------------------------------------------

Since
\[
 a_{\rm TK}=a_m+\rho\psi,
\]
define the distributed thermo-kinetic availability by
\[
 \mathcal A_{\rm TK}(t)
 :=
 \int_0^L a_{\rm TK}(z,t)\,dz .
\]

\begin{proposition}[Distributed relative-energy balance]
\label{prop:global_relative_balance}
For every sufficiently regular solution,
\begin{equation}
\boxed{
\begin{aligned}
 \frac{d}{dt}\mathcal A_{\rm TK}
 ={}&
 Q_m^{\rm rel}(0)-Q_m^{\rm rel}(L)
 +
 Q_{tc}^{\rm rel}(0)-Q_{tc}^{\rm rel}(L)
 \\
 &+
 \int_0^L
 \Big[
 \delta T\,\delta f_j
 +R_r+R_d+R_q
 +R_{\rm ref}+R_{\rm mt}
 \Big]\,dz .
\end{aligned}
}
\label{eq:global_relative_balance}
\end{equation}
\end{proposition}

\begin{proof}
Integrating
\eqref{eq:mechanical_relative_identity} and
\eqref{eq:exact_tc_identity} over $[0,L]$ and adding the resulting
identities gives \eqref{eq:global_relative_balance}.
\end{proof}

At the outlet, the physical conditions
\[
 J(L,t)=J_d(L)=0,
 \qquad
 q(L,t)=q_d(L)=0
\]
reduce the thermochemical relative flux to
\[
 Q_{tc}^{\rm rel}(L,t)
 =
 m(L,t)\psi(L,t).
\]
Under the convexity conditions of
Proposition~\ref{prop:Hamiltonian_convexity},
$\psi\ge0$ in the admissible convex neighborhood of the reference
state. Hence positive outlet flow implies
\[
 Q_{tc}^{\rm rel}(L,t)\ge0.
\]

No analogous unconditional sign is available at the inlet because
the prescribed boundary data are species and total-energy fluxes
rather than the thermochemical state itself. The inlet contribution
is therefore treated together with the bulk thermochemical terms in
the local coercivity analysis of Section~\ref{sec:control}.

For comparison, if the reference thermochemical state is spatially
uniform and the irreversible fluxes vanish, then
$R_{\rm ref}=R_{\rm mt}=0$ and
\[
 R_r=-T_d\,\delta r\,\delta X_r,
 \qquad
 R_d=-T_d\,\delta J\,\delta X_d,
 \qquad
 R_q
 =
 -kT_d
 \left|
 \partial_z\left(\frac{T}{T_d}\right)
 \right|^2
 \left(\frac{T_d}{T}\right)^2 .
\]
Thus the standard incremental dissipative structure is recovered.
For the spatially varying stationary profile considered here, the
additional terms collected in \eqref{eq:tc_remainders} must instead
be retained.

Equation~\eqref{eq:global_relative_balance} is the starting point for
the control design. The jacket enters through the shifted supply
$\delta T\,\delta f_j$, whereas the mechanical subsystem remains
autonomous. The feedback must therefore supply sufficient
thermochemical dissipation to dominate the locally non-dissipative
profile-dependent terms.

%===========================================================
\section{Bounded Monotone Jacket Feedback and Stability Estimates}
\label{sec:control}
%===========================================================

The relative-energy balance of
Section~\ref{sec:relative_balance} shows that the distributed jacket
enters the thermo-kinetic availability through the shifted thermal
supply $\delta T\,\delta f_j$, with $f_j=q_j/T$.
We first construct a bounded monotone feedback that makes this
contribution nonpositive while leaving the intrinsic reaction,
diffusion, conduction, and convection mechanisms unchanged.
The resulting local dissipation is then used in the coercivity and
stability estimates below.

%-----------------------------------------------------------
\subsection{Bounded monotone jacket feedback}
\label{subsec:bounded_feedback}
%-----------------------------------------------------------

The jacket heat exchange is
\[
 q_j=\kappa_j(T_j-T),
 \qquad
 q_{j,d}=\kappa_j(T_{j,d}-T_d),
 \qquad
 \kappa_j>0 .
\]

Consider
\begin{equation}
 \boxed{
 T_j(z,t)
 =
 T_{j,d}(z)
 -
 \varphi\!\left(z,\delta T(z,t)\right).
 }
 \label{eq:general_monotone_feedback}
\end{equation}
For each $z$, assume that $s\mapsto\varphi(z,s)$ is nondecreasing and
\[
 \varphi(z,0)=0 .
\]
It follows immediately that
\[
 s\,\varphi(z,s)\ge0
\]
for every admissible $s$.

\begin{proposition}[Exact shifted jacket dissipation]
\label{prop:jacket_dissipation}
Assume
$T>0$, $T_d>0$, and $T_{j,d}>0$.
Under \eqref{eq:general_monotone_feedback},
\[
 \delta f_j
 =
 -\frac{\kappa_j}{T}
 \left[
 \frac{T_{j,d}}{T_d}\delta T
 +
 \varphi(z,\delta T)
 \right],
\]
and therefore
\begin{equation}
 \boxed{
 \delta T\,\delta f_j
 =
 -\frac{\kappa_j}{T}
 \left[
 \frac{T_{j,d}}{T_d}(\delta T)^2
 +
 \delta T\,\varphi(z,\delta T)
 \right]
 \le0 .
 }
 \label{eq:exact_shifted_jacket_dissipation}
\end{equation}
\end{proposition}

\begin{proof}
Using the jacket law and
\eqref{eq:general_monotone_feedback},
\[
 f_j
 =
 \frac{\kappa_j}{T}
 \left[
 T_{j,d}-\varphi(z,\delta T)-T
 \right],
 \qquad
 f_{j,d}
 =
 \frac{\kappa_j}{T_d}(T_{j,d}-T_d).
\]
Subtracting the two expressions yields the stated formula for
$\delta f_j$. Multiplication by $\delta T$, together with
$T_{j,d}/T_d>0$ and
$\delta T\,\varphi(z,\delta T)\ge0$, gives
\eqref{eq:exact_shifted_jacket_dissipation}.
\end{proof}

Even with $\varphi\equiv0$, the stationary jacket contributes the
shifted dissipation
\[
 \delta T\,\delta f_j
 =
 -\kappa_j
 \frac{T_{j,d}}{TT_d}
 (\delta T)^2
 \le0 .
\]
Thus the feedback adds nonnegative thermal damping to the dissipation
already induced by the stationary jacket.

A bounded realization is
\begin{equation}
 \boxed{
 \varphi(z,s)
 =
 \Delta T_j
 \tanh
 \left(
 \frac{\gamma(z)s}{\Delta T_j}
 \right),
 \qquad
 \gamma(z)\ge0,
 \quad
 \Delta T_j>0 .
 }
 \label{eq:tanh_feedback}
\end{equation}
This map is nondecreasing, satisfies $\varphi(z,0)=0$, and obeys
$|\varphi(z,s)|<\Delta T_j$. Hence
\begin{equation}
 T_{j,d}(z)-\Delta T_j
 <
 T_j(z,t)
 <
 T_{j,d}(z)+\Delta T_j .
 \label{eq:jacket_bound}
\end{equation}
If
\[
 \Delta T_j
 <
 \inf_{z\in[0,L]}T_{j,d}(z),
\]
then $T_j$ remains positive.

Finally,
\[
 \left.
 \frac{\partial\varphi}{\partial s}
 \right|_{s=0}
 =
 \gamma(z).
\]
Thus $\gamma(z)$ sets the local feedback slope about the stationary
profile, whereas $\Delta T_j$ sets the saturation amplitude of the
jacket-temperature correction.

%-----------------------------------------------------------
\subsection{Local thermochemical coercivity and gain condition}
\label{subsec:coercivity_gain}
%-----------------------------------------------------------

We now characterize the quadratic thermochemical dissipation about the
stationary profile. In this subsection the mechanical perturbations are
set to zero, $\delta\rho=\delta m=0$; their coupling with the
thermochemical variables is treated separately in the stability
analysis below.

Let
\[
 \vartheta
 =
 \begin{bmatrix}
 \delta T\\
 \delta\alpha
 \end{bmatrix},
 \qquad
 \xi
 =
 \begin{bmatrix}
 \delta T\\
 \delta\alpha\\
 \delta T_z\\
 \delta\alpha_z
 \end{bmatrix}.
\]
The quadratic thermochemical storage is
\[
 \mathcal A_{tc}^{(2)}
 =
 \frac12
 \int_0^L
 \rho_d
 \left[
 \frac{c_{p,d}}{T_d}(\delta T)^2
 +
 A_{\alpha,d}(\delta\alpha)^2
 \right]dz,
\]
where
\[
 A_{\alpha,d}
 =
 \left.
 \frac{\partial A_r}{\partial\alpha}
 \right|_{T,d}
 >0,
 \qquad
 A_{\alpha,d}
 =
 \frac{R_sT_d}{\alpha_d(1-\alpha_d)}
\]
for the ideal-mixture closure. Throughout this subsection,
$\operatorname{sym}(M):=(M+M^\top)/2$, and all constitutive
coefficients are evaluated along the stationary profile unless stated
otherwise.

For the reaction contribution, define
\[
 a_r
 :=
 \begin{bmatrix}
 r_T\\ r_\alpha
 \end{bmatrix}_d,
 \qquad
 b_r
 :=
 \begin{bmatrix}
 X_{r,T}\\ X_{r,\alpha}
 \end{bmatrix}_d,
 \qquad
 e_T
 :=
 \begin{bmatrix}
 1\\0
 \end{bmatrix}.
\]
The quadratic part of $R_r$ in \eqref{eq:tc_remainders} is
\[
 R_r^{(2)}
 =
 -\vartheta^\top D_r\vartheta,
 \qquad
 D_r
 =
 T_d\operatorname{sym}(a_rb_r^\top)
 -
 r_d\operatorname{sym}(e_Tb_r^\top).
\]

For diffusion, write
\[
 J=L X_d,
 \qquad
 L:=M_dT
 =
 \frac{D\rho\alpha(1-\alpha)}{R_s}.
\]
On the mechanical reference manifold,
\[
 \delta X_d
 =
 b_d^\top\vartheta+B_d\vartheta_z
 =
 g_d^\top\xi,
\]
where
\[
 B_d
 =
 -\frac1{T_d}
 \begin{bmatrix}
 A_{T,d} & A_{\alpha,d}
 \end{bmatrix},
 \qquad
 b_d
 =
 -\frac1{T_d}
 \begin{bmatrix}
 A_{T,d}'+X_{d,d}\\
 A_{\alpha,d}'
 \end{bmatrix},
 \qquad
 g_d
 =
 \begin{bmatrix}
 b_d\\ B_d^\top
 \end{bmatrix}.
\]
Here
$A_T=(\partial A_r/\partial T)_\alpha$.
With
\[
 L_d:=L(\rho_d,\alpha_d),
 \qquad
 L_{\alpha,d}
 :=
 \left.
 \frac{\partial L}{\partial\alpha}
 \right|_d,
 \qquad
 c_d
 =
 \begin{bmatrix}
 J_d\\
 -T_dX_{d,d}L_{\alpha,d}\\
 0\\
 0
 \end{bmatrix},
\]
the quadratic diffusion contribution is
\[
 R_d^{(2)}
 =
 -\xi^\top D_d\xi,
 \qquad
 D_d
 =
 T_dL_dg_dg_d^\top
 -
 \operatorname{sym}(c_dg_d^\top).
\]

For the constant-conductivity Fourier law, let
\[
 g_q
 =
 \begin{bmatrix}
 -T_d'/T_d^2\\
 0\\
 1/T_d\\
 0
 \end{bmatrix},
 \qquad
 c_q
 =
 \begin{bmatrix}
 kT_d'/T_d\\
 0\\
 0\\
 0
 \end{bmatrix}.
\]
The quadratic part of \eqref{eq:Rq_ratio} is
\[
 R_q^{(2)}
 =
 -\xi^\top D_q\xi,
 \qquad
 D_q
 =
 kT_dg_qg_q^\top
 -
 \operatorname{sym}(c_qg_q^\top).
\]

The spatially varying reference also contributes at second order.
The local coordinate transformation satisfies
\[
 \delta\zeta
 =
 S_d\vartheta+\mathcal O(|\vartheta|^2),
 \qquad
 S_d
 =
 \begin{bmatrix}
 c_{p,d}/T_d & -A_{T,d}\\
 0 & 1
 \end{bmatrix}.
\]
Since
\[
 r_w
 =
 \frac12
 \nabla_\zeta^3u_{tc}(\zeta_d)
 [\delta\zeta,\delta\zeta]
 +
 \mathcal O(|\delta\zeta|^3),
\]
symmetry of the third derivative gives
\[
 R_{\rm ref}^{(2)}
 =
 \frac{m_d}{2}
 \delta\zeta^\top H_d'\delta\zeta
 =
 -\xi^\top D_{\rm ref}\xi,
\]
with
\[
 D_{\rm ref}
 =
 -\frac{m_d}{2}
 \begin{bmatrix}
 S_d^\top H_d'S_d & 0_{2\times2}\\
 0_{2\times2} & 0_{2\times2}
 \end{bmatrix}.
\]

Finally, linearization of the bounded feedback
\eqref{eq:tanh_feedback} gives
\[
 \delta T\,\delta f_j
 =
 -K_{j,d}(\delta T)^2
 +
 \mathcal O(|\delta T|^3),
\]
where
\[
 K_{j,d}
 =
 \kappa_j
 \left(
 \frac{T_{j,d}}{T_d^2}
 +
 \frac{\gamma}{T_d}
 \right).
\]
The first term,
$\kappa_jT_{j,d}/T_d^2$, is the natural shifted jacket damping,
whereas $\kappa_j\gamma/T_d$ is the additional local damping supplied
by the feedback.

Embedding
\[
 \widehat D_r
 =
 \begin{bmatrix}
 D_r & 0_{2\times2}\\
 0_{2\times2} & 0_{2\times2}
 \end{bmatrix},
 \qquad
 E_T=
 \begin{bmatrix}
 1&0&0&0
 \end{bmatrix}^{\top},
\]
the complete thermochemical bulk matrix is
\begin{equation}
 \boxed{
 \mathbb D_\gamma
 =
 \widehat D_r
 +
 D_d
 +
 D_q
 +
 D_{\rm ref}
 +
 K_{j,d}E_TE_T^\top .
 }
 \label{eq:Dgamma}
\end{equation}

The bulk matrix alone does not characterize stability of the
distributed system because the relative-energy identity also contains
boundary quadratic forms. Let
$\mathcal B_{\rm in}^{(2)}[\vartheta(0)]$ denote the thermochemical
second variation of the inlet relative flux after imposing the
linearized physical inlet constraints. At the outlet,
\[
 \mathcal B_{\rm out}^{(2)}[\vartheta(L)]
 =
 \frac{m_d}{2}
 \vartheta(L)^\top G_d(L)\vartheta(L),
 \qquad
 G_d
 =
 \begin{bmatrix}
 c_{p,d}/T_d & 0\\
 0 & A_{\alpha,d}
 \end{bmatrix}.
\]
Define the complete quadratic dissipation functional
\begin{equation}
\begin{aligned}
 \mathfrak D_\gamma[\vartheta]
 :={}&
 \int_0^L
 \xi^\top\mathbb D_\gamma\xi\,dz
 -
 \mathcal B_{\rm in}^{(2)}[\vartheta(0)]
 +
 \mathcal B_{\rm out}^{(2)}[\vartheta(L)] .
\end{aligned}
\label{eq:global_dissipation_functional}
\end{equation}

\begin{assumption}[Local thermochemical coercivity]
\label{ass:tc_coercivity}
There exists $\nu_{tc}>0$ such that
\begin{equation}
 \boxed{
 \mathfrak D_\gamma[\vartheta]
 \ge
 \nu_{tc}
 \|\vartheta\|_{H^1(0,L)}^2
 }
 \label{eq:tc_coercivity}
\end{equation}
for every admissible thermochemical perturbation satisfying the
linearized physical boundary constraints.
\end{assumption}

Assumption~\ref{ass:tc_coercivity} is the continuous bulk--boundary
condition used in the stability results. For gain selection, a simpler
pointwise condition can be obtained from the bulk matrix.

Let $\mathbb D_0$ denote \eqref{eq:Dgamma} with $\gamma=0$ and write
\[
 \mathbb D_0
 =
 \begin{bmatrix}
 d_{11}^0 & d_{12}^\top\\
 d_{12} & D_{22}
 \end{bmatrix}.
\]
Since the feedback acts only on the temperature coordinate,
\[
 \mathbb D_\gamma
 =
 \begin{bmatrix}
 d_{11}^0+\kappa_j\gamma/T_d & d_{12}^\top\\
 d_{12} & D_{22}
 \end{bmatrix}.
\]

\begin{proposition}[Pointwise sufficient condition for positive definiteness of the bulk matrix]
\label{prop:pointwise_gain}
Assume
$D_{22}(z)\succ0$ for all $z\in[0,L]$.
Then
$\mathbb D_\gamma(z)\succ0$ whenever
\begin{equation}
 \boxed{
 \gamma(z)
 >
 \gamma^\star(z)
 :=
 \frac{T_d(z)}{\kappa_j}
 \left[
 d_{12}(z)^\top
 D_{22}(z)^{-1}
 d_{12}(z)
 -
 d_{11}^0(z)
 \right].
 }
 \label{eq:gamma_star}
\end{equation}
Together with the monotonicity requirement, the admissible gain must
therefore satisfy
\[
 \gamma(z)\ge0,
 \qquad
 \gamma(z)>\gamma^\star(z).
\]
\end{proposition}

\begin{proof}
The Schur-complement criterion gives
$\mathbb D_\gamma\succ0$ if and only if
$D_{22}\succ0$ and
\[
 d_{11}^0
 +
 \frac{\kappa_j}{T_d}\gamma
 -
 d_{12}^\top D_{22}^{-1}d_{12}
 >0,
\]
which is equivalent to \eqref{eq:gamma_star}.
\end{proof}

Condition~\eqref{eq:gamma_star} guarantees pointwise positivity of the
bulk matrix only. It does not, by itself, imply
Assumption~\ref{ass:tc_coercivity}, because the latter also contains
the physical boundary contributions. The pointwise condition is
therefore used as a convenient conservative gain-design criterion,
whereas the full bulk--boundary coercivity remains the condition
required by the continuous stability analysis. Generalized-eigenvalue
tests after structure-preserving discretization are used only to
assess the corresponding semi-discrete quadratic functional and are
not interpreted as a proof of continuous coercivity.

%-----------------------------------------------------------
\subsection{Conditional nonlinear relative-energy estimate}
\label{subsec:nonlinear_estimate}
%-----------------------------------------------------------

The exact jacket dissipation
\eqref{eq:exact_shifted_jacket_dissipation} is nonlinear, whereas
Assumption~\ref{ass:tc_coercivity} concerns the quadratic
thermochemical dissipation about the stationary profile. We now relate
the two locally.

Define
\[
 \mathcal A_{tc}(t)
 :=
 \int_0^L\rho\psi\,dz,
 \qquad
 \mathcal B_m(t)
 :=
 Q_m^{\rm rel}(L,t)-Q_m^{\rm rel}(0,t)\ge0.
\]
By Proposition~\ref{prop:mechanical_boundary},
$\dot{\mathcal A}_m=-\mathcal B_m$.

\begin{theorem}[Conditional nonlinear relative-energy estimate]
\label{thm:nonlinear_estimate}
Assume that the constitutive functions are $C^3$ in a neighborhood of
the stationary profile, that the profile and its first two spatial
derivatives are bounded, and that
Assumption~\ref{ass:tc_coercivity} holds. Suppose further that a
sufficiently regular closed-loop classical solution exists and remains
in a sufficiently small compact admissible neighborhood satisfying
\[
 \rho\ge\rho_{\min}>0,
 \qquad
 T\ge T_{\min}>0,
 \qquad
 0<\alpha_{\min}\le\alpha\le\alpha_{\max}<1,
 \qquad
 0<v<c(\rho).
\]
Then there exist constants
$c_{tc}>0$, $C_m>0$, and $C_b\ge0$ such that
\begin{equation}
 \dot{\mathcal A}_{tc}
 \le
 -c_{tc}\|\vartheta\|_{H^1(0,L)}^2
 +
 C_m\mathcal A_m
 +
 C_b\mathcal B_m .
 \label{eq:nonlinear_Atc_estimate}
\end{equation}
Consequently, for any $\Gamma_m>C_b$, the functional
\[
 W:=\mathcal A_{tc}+\Gamma_m\mathcal A_m
\]
satisfies
\begin{equation}
 \boxed{
 \dot W
 \le
 -c_{tc}\|\vartheta\|_{H^1(0,L)}^2
 +
 C_m\mathcal A_m .
 }
 \label{eq:nonlinear_composite_estimate}
\end{equation}
\end{theorem}

\begin{proof}
Expand the exact relative-energy balance
\eqref{eq:global_relative_balance} about the stationary profile.
Assumption~\ref{ass:tc_coercivity} bounds the purely thermochemical
quadratic part above by
$-\nu_{tc}\|\vartheta\|_{H^1}^2$.

The higher-order thermochemical terms satisfy locally
\[
 |\mathcal R_{\rm nl}|
 \le
 C\|\vartheta\|_{L^\infty}
 \|\vartheta\|_{H^1}^2.
\]
Since $H^1(0,L)\hookrightarrow L^\infty(0,L)$, the neighborhood can be
chosen sufficiently small that these terms are absorbed into the
coercive quadratic contribution.

The terms containing the mechanical perturbations, including
$R_{\rm mt}$ and the mechanical dependence of the constitutive laws,
are bilinear to leading order. Local equivalence of the mechanical
relative energy with the $L^2$ norm of
$(\delta\rho,\delta m)$, together with Young's inequality, gives
\[
 |\mathcal C_{\rm bulk}|
 \le
 \varepsilon\|\vartheta\|_{H^1}^2
 +
 C_\varepsilon\mathcal A_m .
\]

At the inlet, strict subsonicity and continuity imply, after reducing
the neighborhood if necessary,
\[
 -Q_m^{\rm rel}(0)
 \ge
 c_0|\delta\rho(0)|^2,
 \qquad
 Q_m^{\rm rel}(L)
 \ge
 c_L|\delta v(L)|^2
\]
for some $c_0,c_L>0$. Hence the mechanical boundary dissipation
controls the mechanical traces entering the thermochemical coupling,
and
\[
 |\mathcal C_{\rm in}|
 \le
 \varepsilon|\vartheta(0)|^2
 +
 C_\varepsilon\mathcal B_m .
\]
The thermochemical trace term is absorbed by
Assumption~\ref{ass:tc_coercivity} and the trace inequality.
Choosing the neighborhood and $\varepsilon$ sufficiently small yields
\eqref{eq:nonlinear_Atc_estimate}.

Finally,
\[
 \dot W
 \le
 -c_{tc}\|\vartheta\|_{H^1}^2
 +
 C_m\mathcal A_m
 -
 (\Gamma_m-C_b)\mathcal B_m ,
\]
and \eqref{eq:nonlinear_composite_estimate} follows from
$\Gamma_m>C_b$.
\end{proof}

Theorem~\ref{thm:nonlinear_estimate} is an a-priori estimate for
sufficiently regular solutions and does not establish nonlinear
well-posedness.

\begin{corollary}[Invariant mechanical manifold]
\label{cor:mechanical_invariant_manifold}
Suppose
\[
 \rho(z,0)=\rho_d,
 \qquad
 m(z,0)=m_d,
 \qquad
 z\in[0,L].
\]
For every classical solution for which uniqueness holds, autonomy of
the mechanical subsystem
\eqref{eq:model_mass}--\eqref{eq:model_momentum} implies
\[
 \rho(z,t)=\rho_d,
 \qquad
 m(z,t)=m_d
\]
throughout the interval of existence.

Under the assumptions of
Theorem~\ref{thm:nonlinear_estimate}, there exists
$\lambda_{tc}>0$ such that
\begin{equation}
 \boxed{
 \mathcal A_{tc}(t)
 \le
 e^{-\lambda_{tc}t}\mathcal A_{tc}(0).
 }
 \label{eq:nonlinear_tc_exponential}
\end{equation}
\end{corollary}

\begin{proof}
The constant state $(\rho_d,m_d)$ satisfies the autonomous mechanical
equations and the same mechanical boundary conditions. Uniqueness
therefore implies invariance, so that
$\mathcal A_m=\mathcal B_m=0$.
Equation~\eqref{eq:nonlinear_Atc_estimate} reduces to
\[
 \dot{\mathcal A}_{tc}
 \le
 -c_{tc}\|\vartheta\|_{H^1}^2.
\]
Local strict convexity gives
\[
 \mathcal A_{tc}
 \le
 C_A\|\vartheta\|_{L^2}^2
 \le
 C_A\|\vartheta\|_{H^1}^2,
\]
and hence
$\dot{\mathcal A}_{tc}\le
-(c_{tc}/C_A)\mathcal A_{tc}$.
Integration yields \eqref{eq:nonlinear_tc_exponential}.
\end{proof}

This is a conditional thermochemical exponential-decay result on a invariant mechanical manifold; no general nonlinear exponential
stability of the full PDE is claimed.

%-----------------------------------------------------------
\subsection{Exponential estimate for the linearized closed-loop PDE}
\label{subsec:linearized_estimate}
%-----------------------------------------------------------

We finally derive an energy estimate for the complete linearized
closed-loop system. Let
\[
 \varrho:=\delta\rho,
 \qquad
 u:=\delta v,
 \qquad
 c_s^2:=p'(\rho_d).
\]
The linearized mechanical equations are
\[
 \varrho_t+v_d\varrho_z+\rho_du_z=0,
 \qquad
 u_t+v_du_z+\frac{c_s^2}{\rho_d}\varrho_z=0.
\]

Introduce the Riemann variables
\[
 w_1=u+\frac{c_s}{\rho_d}\varrho,
 \qquad
 w_2=u-\frac{c_s}{\rho_d}\varrho,
\]
which satisfy
\[
 (w_1)_t+(v_d+c_s)(w_1)_z=0,
 \qquad
 (w_2)_t+(v_d-c_s)(w_2)_z=0.
\]
The linearized boundary conditions become
\[
 w_1(0)=-\kappa_bw_2(0),
 \qquad
 \kappa_b=\frac{c_s-v_d}{c_s+v_d},
 \qquad
 w_2(L)=w_1(L).
\]

For $\mu>0$, define
\begin{equation}
 E_m^\mu
 =
 \frac{\rho_d}{4}
 \int_0^L
 \left[
 e^{-\mu z}w_1^2
 +
 e^{\mu z}w_2^2
 \right]dz .
 \label{eq:weighted_mechanical_energy}
\end{equation}

\begin{lemma}[Linear mechanical decay]
\label{lem:linear_mechanical_decay}
If
\begin{equation}
 0<\mu<
 \frac{1}{2L}
 \ln\left(
 \frac{c_s+v_d}{c_s-v_d}
 \right),
 \label{eq:mu_condition}
\end{equation}
then
\begin{equation}
 \dot E_m^\mu
 \le
 -\lambda_mE_m^\mu
 -
 \beta_0w_2^2(0)
 -
 \beta_Lw_1^2(L),
 \label{eq:mechanical_weighted_decay}
\end{equation}
where
\[
 \lambda_m=\mu(c_s-v_d)>0,
\]
\[
 \beta_0
 =
 \frac{\rho_dv_d(c_s-v_d)}
 {2(c_s+v_d)}
 >0,
\]
and
\[
 \beta_L
 =
 \frac{\rho_d}{4}
 \left[
 (c_s+v_d)e^{-\mu L}
 -
 (c_s-v_d)e^{\mu L}
 \right]
 >0.
\]
\end{lemma}

\begin{proof}
Differentiating \eqref{eq:weighted_mechanical_energy} along the two
transport equations and integrating by parts gives
\begin{align*}
 \dot E_m^\mu
 ={}&
 -\frac{\rho_d\mu}{4}
 \int_0^L
 \left[
 (c_s+v_d)e^{-\mu z}w_1^2
 +
 (c_s-v_d)e^{\mu z}w_2^2
 \right]dz
 \\
 &+
 \frac{\rho_d}{4}
 \left[
 (c_s+v_d)w_1^2(0)
 -
 (c_s-v_d)w_2^2(0)
 \right]
 \\
 &-
 \frac{\rho_d}{4}
 \left[
 (c_s+v_d)e^{-\mu L}w_1^2(L)
 -
 (c_s-v_d)e^{\mu L}w_2^2(L)
 \right].
\end{align*}
The inlet relation gives the coefficient $-\beta_0$, whereas the
outlet relation gives $-\beta_L$.
Condition~\eqref{eq:mu_condition} guarantees $\beta_L>0$, and the
distributed term is bounded above by
$-\mu(c_s-v_d)E_m^\mu$.
\end{proof}

The thermochemical linearization is governed by the quadratic
bulk--boundary functional of
Assumption~\ref{ass:tc_coercivity}, with the mechanical variables
entering as an inhomogeneous cascade input.

\begin{lemma}[Linear thermochemical estimate]
\label{lem:linear_tc_estimate}
Under Assumption~\ref{ass:tc_coercivity}, there exist constants
$\lambda_{tc}^{\rm lin}>0$,
$C_m^{\rm lin}\ge0$, and
$C_b^{\rm lin}\ge0$ such that every sufficiently regular solution
satisfies
\begin{equation}
 \dot{\mathcal A}_{tc}^{(2)}
 \le
 -\lambda_{tc}^{\rm lin}\mathcal A_{tc}^{(2)}
 +
 C_m^{\rm lin}E_m^\mu
 +
 C_b^{\rm lin}w_2^2(0).
 \label{eq:linear_tc_energy_estimate}
\end{equation}
\end{lemma}

\begin{proof}
The second variation of
\eqref{eq:exact_tc_identity} yields the quadratic thermochemical
bulk--boundary functional of
Assumption~\ref{ass:tc_coercivity}.
The mechanical coupling terms satisfy, for every $\varepsilon>0$,
\[
 |\mathcal C_{\rm bulk}|
 \le
 \varepsilon\|\vartheta\|_{H^1}^2
 +
 C_\varepsilon E_m^\mu ,
\]
while the inlet coupling satisfies
\[
 |\mathcal C_{\rm in}|
 \le
 \varepsilon|\vartheta(0)|^2
 +
 C_\varepsilon w_2^2(0).
\]
The terms proportional to $\varepsilon$ are absorbed by
Assumption~\ref{ass:tc_coercivity} and the trace inequality.
Local equivalence of
$\mathcal A_{tc}^{(2)}$ with
$\|\vartheta\|_{L^2}^2$
then gives \eqref{eq:linear_tc_energy_estimate}.
\end{proof}

\begin{theorem}[Exponential energy estimate for the linearized closed loop]
\label{thm:linearized_exponential}
Suppose Assumption~\ref{ass:tc_coercivity} holds and choose $\mu$
according to \eqref{eq:mu_condition}. Define
\[
 V_{\rm lin}
 :=
 \mathcal A_{tc}^{(2)}
 +
 \kappa_mE_m^\mu .
\]
If
\begin{equation}
 \boxed{
 \kappa_m
 >
 \max
 \left\{
 \frac{C_m^{\rm lin}}{\lambda_m},
 \frac{C_b^{\rm lin}}{\beta_0}
 \right\},
 }
 \label{eq:kappa_m_condition}
\end{equation}
then there exists $\lambda_V>0$ such that every sufficiently regular
solution of the linearized closed-loop initial-boundary-value problem
satisfies
\begin{equation}
 \boxed{
 \dot V_{\rm lin}\le-\lambda_VV_{\rm lin},
 \qquad
 V_{\rm lin}(t)
 \le
 e^{-\lambda_Vt}V_{\rm lin}(0).
 }
 \label{eq:V_linear_exponential}
\end{equation}
\end{theorem}

\begin{proof}
Combining
\eqref{eq:mechanical_weighted_decay} and
\eqref{eq:linear_tc_energy_estimate} yields
\begin{align*}
 \dot V_{\rm lin}
 \le{}&
 -\lambda_{tc}^{\rm lin}\mathcal A_{tc}^{(2)}
 -
 \left(
 \kappa_m\lambda_m-C_m^{\rm lin}
 \right)E_m^\mu
 \\
 &-
 \left(
 \kappa_m\beta_0-C_b^{\rm lin}
 \right)w_2^2(0)
 -
 \kappa_m\beta_Lw_1^2(L).
\end{align*}
Condition~\eqref{eq:kappa_m_condition} makes all coefficients on the
right-hand side strictly positive. Since $V_{\rm lin}$ is equivalent
to
$\mathcal A_{tc}^{(2)}+E_m^\mu$, the estimate
\eqref{eq:V_linear_exponential} follows.
\end{proof}

Theorem~\ref{thm:linearized_exponential} is an exponential energy
estimate for sufficiently regular solutions of the linearized
closed-loop initial-boundary-value problem. It is not a semigroup
exponential-stability statement, because generation and well-posedness
of the associated linear operator are not established here.

%===========================================================

%===========================================================
\section{Numerical Verification}
\label{sec:numerical}
%===========================================================

The numerical study examines three aspects of the proposed design:
the pointwise bulk-positivity gain criterion, the transient improvement
produced by the bounded jacket feedback, and preservation of the
thermodynamic structure in the structure-preserving discretization.

The open-loop reference case uses the physical reactor with the jacket
fixed at its target value,
\[
 T_j(z,t)=T_{j,d}=370~{\rm K}.
\]
For the selected benchmark this fixed-jacket system is numerically
convergent. The purpose of the feedback is therefore to increase the
structured thermal dissipation and improve the transient rather than
to stabilize an artificially unstable operating point.

%-----------------------------------------------------------
\subsection{Numerical setting and stationary-profile validation}
\label{subsec:numerical_setting}
%-----------------------------------------------------------

The simulations use the structure-preserving
Legendre--Gauss--Lobatto PFEM discretization of the BC-IPHS model.
A single polynomial element of degree $p=60$, corresponding to
$61$ spatial nodes, is used. The numerical and control settings that
are not already specified in Section~\ref{sec:model} are summarized in
Table~\ref{tab:numerical_parameters}.

\begin{table}[t]
\centering
\caption{Numerical and control settings.}
\label{tab:numerical_parameters}
\small
\setlength{\tabcolsep}{4pt}
\begin{tabular}{@{}lcc@{}}
\hline
Quantity & Symbol & Value \\
\hline
Initial-profile jacket, large test
& $T_{j,0}$ & $350~{\rm K}$ \\
Initial-profile jacket, local test
& $T_{j,0}$ & $365~{\rm K}$ \\
Target jacket temperature
& $T_{j,d}$ & $370~{\rm K}$ \\
Feedback slope
& $\gamma$ & $4.3$ \\
Feedback excursion bound
& $\Delta T_j$ & $10~{\rm K}$ \\
Simulation horizon
& $t_f$ & $160~{\rm s}$ \\
Output interval
& $\Delta t$ & $0.5~{\rm s}$ \\
ODE solver & --- & \texttt{ode15s} \\
Relative tolerance & --- & $2\times10^{-9}$ \\
\hline
\end{tabular}
\end{table}

Stationary profiles corresponding to fixed jacket temperatures
$350$, $365$, and $370$ K are computed on the same $p=60$ PFEM
space before the transient simulations. The $370$ K profile defines
the target $(T_d,\alpha_d)$. Comparison with an independently
validated reference computation gives an outlet-temperature
difference of $3.4\times10^{-8}$ K and a product-mass-fraction
difference of $8.5\times10^{-11}$.

Hereafter,
$\alpha_B:=1-\alpha$ denotes the product mass fraction.

For both profile-regulation experiments, the mechanical state is
initialized at
\[
 \rho(z,0)=\rho_d,
 \qquad
 m(z,0)=m_d.
\]
By Corollary~\ref{cor:mechanical_invariant_manifold}, this is an exact
invariant mechanical manifold of the continuous model. The
corresponding PFEM simulations therefore evolve only the
thermochemical dynamics on that invariant manifold.

%-----------------------------------------------------------
\subsection{Gain selection and semi-discrete eigenvalue analysis}
\label{subsec:numerical_gain}
%-----------------------------------------------------------

We first evaluate the pointwise bulk-positivity condition of
Proposition~\ref{prop:pointwise_gain}. Along the target profile,
\[
 \min_{z\in[0,L]}
 \lambda_{\min}\!\left(D_{22}(z)\right)
 \simeq3.29\times10^{-2}>0,
\]
and the required pointwise gain satisfies
\[
 \max_{z\in[0,L]}\gamma^\star(z)<4.16.
\]
The constant choice $\gamma=4.3$ therefore satisfies the derived
pointwise bulk-positivity condition over the complete reactor.

\begin{figure}[t]
\centering
\includegraphics[width=\columnwidth]
{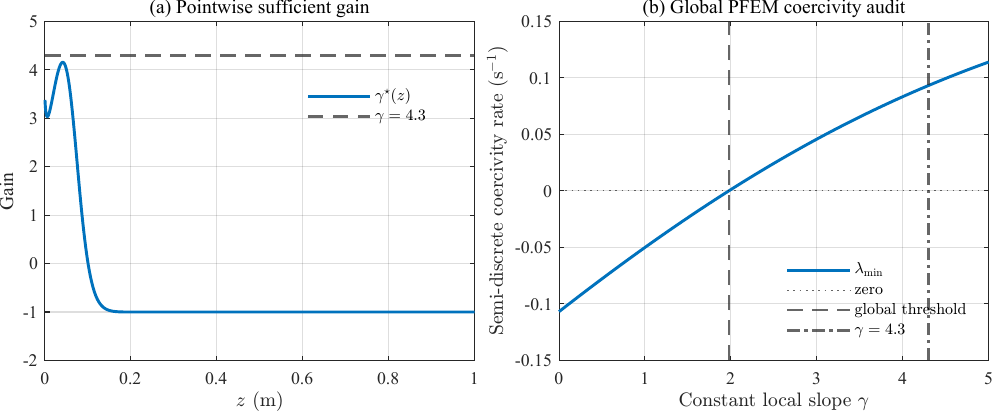}
\caption{Gain selection and semi-discrete eigenvalue analysis.
(a) Pointwise bulk-positivity threshold $\gamma^\star(z)$ and selected
constant slope $\gamma=4.3$.
(b) Lowest generalized dissipation eigenvalue of the $p=60$
semi-discrete PFEM model as a function of the feedback slope.
The zero crossing occurs at approximately $\gamma=1.98$.
The generalized-eigenvalue calculation concerns only the
semi-discrete model and is not a proof of continuous coercivity.}
\label{fig:gain_coercivity}
\end{figure}

Figure~\ref{fig:gain_coercivity}(a) shows the spatial distribution of
$\gamma^\star(z)$. The largest required value occurs close to the
inlet, whereas substantially less additional thermal damping is
required over most of the reactor.

For comparison, the complete semi-discrete bulk--boundary functional
is examined through its lowest generalized eigenvalue.
Figure~\ref{fig:gain_coercivity}(b) shows a zero crossing near
$\gamma=1.98$, considerably below the pointwise threshold.
For this discretization, the difference illustrates the contribution
of the spatial-gradient and boundary terms that are not exploited by
the pointwise Schur-complement test.

This semi-discrete eigenvalue calculation does not replace
Assumption~\ref{ass:tc_coercivity}. The latter remains the continuous
condition used in the stability results.

%-----------------------------------------------------------
\subsection{Distributed profile regulation and disturbance response}
\label{subsec:numerical_regulation}
%-----------------------------------------------------------

All closed-loop simulations use
\[
 T_j(z,t)
 =
 370
 -
 10
 \tanh\!\left(
 \frac{4.3\,[T(z,t)-T_d(z)]}{10}
 \right),
\]
which has the form analyzed in Section~\ref{subsec:bounded_feedback}, rather than a linear command
followed by external saturation.

The first experiment starts from the stationary profile associated
with a $350$ K jacket and regulates the reactor toward the $370$ K
target. Because no quantitative radius is established for the local
nonlinear estimate, this relatively large excursion is regarded as a
performance test rather than as a direct verification of
Theorem~\ref{thm:nonlinear_estimate}.
The $365\rightarrow370$ K experiment provides a substantially smaller
perturbation for comparison.

\begin{figure}[t]
\centering
\includegraphics[width=\columnwidth]
{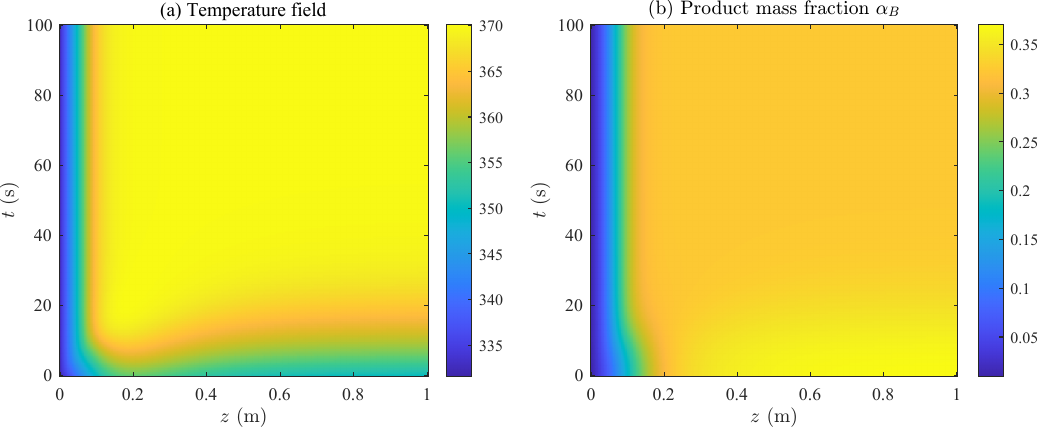}
\caption{Closed-loop distributed-state evolution for the
$350\rightarrow370$ K experiment:
(a) temperature field $T(z,t)$;
(b) product mass fraction $\alpha_B(z,t)$.
Both fields approach the stationary target smoothly.}
\label{fig:distributed_profiles}
\end{figure}

Figure~\ref{fig:distributed_profiles} shows the controlled
temperature and product-composition fields for the
$350\rightarrow370$ K experiment.

At $t_f=160$ s, the maximum spatial errors in the controlled case are
approximately $1.2\times10^{-8}$ K for temperature and
$4.3\times10^{-11}$ for product mass fraction.

\begin{figure}[t]
\centering
\includegraphics[width=\columnwidth]
{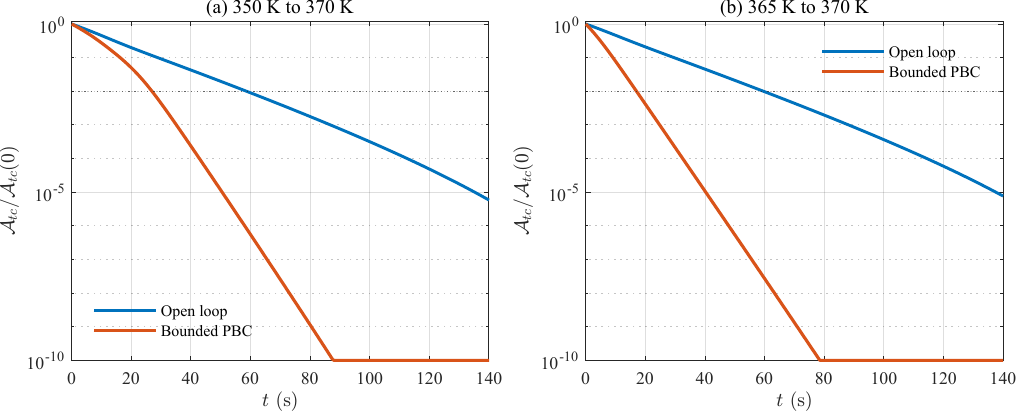}
\caption{Normalized thermochemical relative Hamiltonian:
(a) $350\rightarrow370$ K experiment;
(b) $365\rightarrow370$ K experiment.
The fixed-jacket open loop already exhibits convergence, while the
bounded PBC substantially accelerates the relative-energy decay.}
\label{fig:TKA_decay}
\end{figure}

To compare the transient using the relative storage, define
\[
 \overline{\mathcal A}_{tc}(t)
 :=
 \frac{\mathcal A_{tc}(t)}
 {\mathcal A_{tc}(0)}.
\]
Figure~\ref{fig:TKA_decay} compares the fixed-jacket and
bounded-feedback responses for the $350\rightarrow370$ K and
$365\rightarrow370$ K experiments.

A complementary distributed metric is used to avoid assessing the
transient solely through an integral storage. Let
\[
 e_T^\infty(t)
 :=
 \max_z|T(z,t)-T_d(z)|,
 \qquad
 e_\alpha^\infty(t)
 :=
 \max_z|\alpha_B(z,t)-\alpha_{B,d}(z)|,
\]
and define
\[
 E_\infty(t)
 :=
 \max\left\{
 \frac{e_T^\infty(t)}{e_T^\infty(0)},
 \frac{e_\alpha^\infty(t)}{e_\alpha^\infty(0)}
 \right\}.
\]
For the $350\rightarrow370$ K experiment,
$E_\infty$ first falls below $10^{-2}$ at $131.0$ s in open loop and
$46.5$ s under the bounded PBC. The corresponding relative-energy
and profile-error convergence times are summarized in
Table~\ref{tab:numerical_results}.

To activate the mechanical dynamics, a separate experiment applies a
temporary perturbation to the inlet species and total-energy fluxes:
\[
 \beta(t)=
 \begin{cases}
 1+0.05\sin^2\!\left(\dfrac{\pi(t-10)}{40}\right),
 &10\le t\le50~{\rm s},\\
 1, &\text{otherwise}.
 \end{cases}
\]
All prescribed inlet fluxes are multiplied by the same factor, so
the inlet composition and specific energy are unchanged while the
mass throughput undergoes a smooth temporary $5\%$ increase.
No velocity boundary condition is added; the velocity response follows
from the mass--momentum dynamics through $v=m/\rho$.

\begin{figure}[t]
\centering
\includegraphics[width=\columnwidth]
{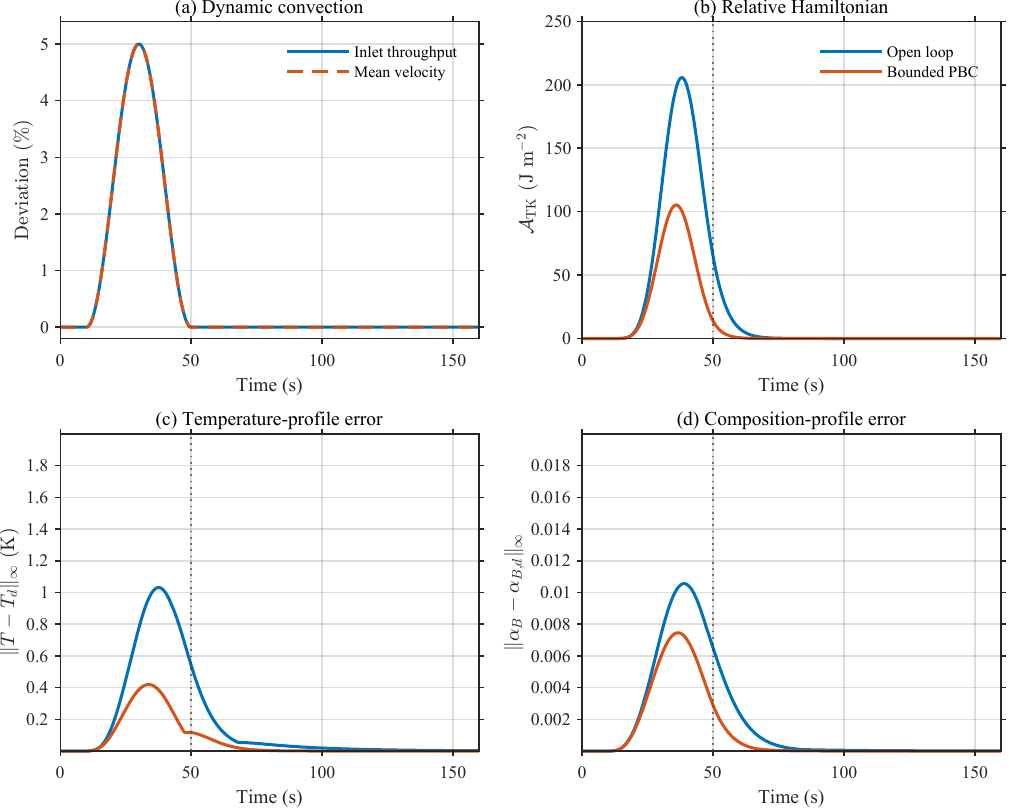}
\caption{Response to a temporary $5\%$ inlet-throughput perturbation:
(a) throughput factor and mean convection velocity;
(b) full thermo-kinetic relative Hamiltonian;
(c) maximum temperature error;
(d) maximum product-composition error.}
\label{fig:throughput_disturbance}
\end{figure}

As shown in Fig.~\ref{fig:throughput_disturbance}, the mean velocity
increases from approximately
$5.00\times10^{-3}$ to $5.25\times10^{-3}$ m\,s$^{-1}$ and is
essentially identical with and without jacket feedback, consistently
with the autonomous mechanical subsystem. The bounded PBC nevertheless
attenuates the thermochemical response: the peak temperature error is
reduced from $1.03$ to $0.42$ K, the peak composition error from
$1.06\times10^{-2}$ to $7.46\times10^{-3}$, and the peak full
relative Hamiltonian from $205.7$ to $105.3$ J\,m$^{-2}$.

During the perturbation, the inlet boundary data differ from those of
the stationary reference. This experiment is therefore not a direct
verification of the fixed-boundary stability estimates; it illustrates how the mechanical perturbation propagates into the
thermochemical variables.

%-----------------------------------------------------------
\subsection{Bounded actuation and thermodynamic diagnostics}
\label{subsec:numerical_structure}
%-----------------------------------------------------------

\begin{figure}[t]
\centering
\includegraphics[width=\columnwidth]
{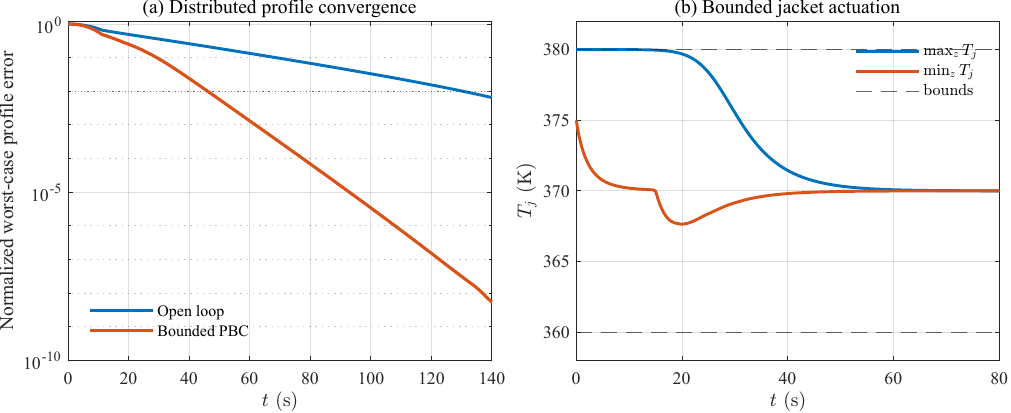}
\caption{Closed-loop performance and actuation for the
$350\rightarrow370$ K experiment.
(a) Normalized worst-case distributed profile error $E_\infty(t)$.
(b) Maximum and minimum jacket temperatures.
The theoretical bounds $360<T_j<380$ K are respected throughout the
transient.}
\label{fig:tracking_actuation}
\end{figure}

Figure~\ref{fig:tracking_actuation} compares the worst-case
distributed profile error with the jacket-temperature extrema during
the $350\rightarrow370$ K experiment.

The realized jacket temperature lies approximately in
$[367.65,380)$ K, while the controller construction guarantees
$360<T_j<380$ K. Thus the actuator bound is enforced directly by the
feedback law rather than by post-processing or command clipping.

The shifted jacket supply of
Proposition~\ref{prop:jacket_dissipation} also remains nonpositive
along the computed trajectory:
\[
 \int_0^L\delta T\,\delta f_j\,dz\le0,
\]
up to machine precision.

Finally, the scaled first- and second-law residuals satisfy
\[
 \max_t r_H(t)=8.2\times10^{-13},
 \qquad
 \max_t r_S(t)=3.6\times10^{-14}.
\]
The integrated reaction, diffusion--heat, and jacket entropy-production
contributions remain nonnegative throughout the controlled transient.
These diagnostics confirm preservation of the discrete first- and
second-law structure of the PFEM simulation.

\begin{table}[t]
\centering
\caption{Closed-loop performance of the bounded PBC.}
\label{tab:numerical_results}
\small
\setlength{\tabcolsep}{4pt}
\begin{tabular}{@{}lcc@{}}
\hline
Metric & Open loop & Bounded PBC \\
\hline
$350\!\rightarrow\!370$ K: $t_{\mathcal A,1\%}$
& $59.0$ s & $27.0$ s \\
$350\!\rightarrow\!370$ K: $t_{E,1\%}$
& $131.0$ s & $46.5$ s \\
$365\!\rightarrow\!370$ K: $t_{\mathcal A,1\%}$
& $60.0$ s & $17.5$ s \\
$350\!\rightarrow\!370$ K: $T_j$ range
& $370$ K & $[367.65,380)$ K \\
\hline
\end{tabular}
\end{table}

Here $t_{\mathcal A,1\%}$ and $t_{E,1\%}$ denote the first times at
which $\overline{\mathcal A}_{tc}$ and $E_\infty$, respectively, fall
below $10^{-2}$.

Overall, the numerical results show that the selected bounded feedback
satisfies the pointwise bulk-positivity gain criterion, substantially
accelerates the thermochemical transient relative to the numerically
convergent fixed-jacket case, and preserves both the actuator bound and
the discrete thermodynamic structure.

%===========================================================
\section{Conclusion}
\label{sec:conclusion}
%===========================================================

This paper developed a relative-energy control approach for a
boundary-controlled irreversible port-Hamiltonian tubular reactor with
dynamic momentum and a spatially varying stationary thermochemical
profile. The thermo-kinetic availability was identified with the
relative Hamiltonian of the BC-IPHS and, for the adopted constitutive
closure, was shown to be positive definite locally and to decompose
into kinetic, compressive, thermal, and mixing contributions. An exact
relative-energy balance was then derived, explicitly retaining the
terms generated by the nonuniform stationary profile.

A bounded monotone jacket feedback was constructed whose shifted
thermal supply is nonpositive exactly. The continuous stability
analysis was formulated in terms of a bulk--boundary thermochemical
coercivity condition, while a Schur-complement argument provided a
pointwise sufficient condition for positive definiteness of the bulk matrix and a practical
gain-selection rule. For sufficiently regular classical solutions
remaining in a small admissible neighborhood, a conditional nonlinear
relative-energy estimate was obtained. Thermochemical exponential
decay follows on the invariant mechanical manifold, and an exponential
energy estimate was established for sufficiently regular solutions of
the complete linearized closed-loop initial-boundary-value problem.

The structure-preserving PFEM simulations support the proposed gain
selection and show substantial acceleration of the thermochemical
transient while respecting the prescribed jacket-temperature bound and
preserving the discrete thermodynamic identities. The selected
fixed-jacket benchmark is already numerically convergent, so the
controller acts to increase thermal dissipation rather than
to stabilize an artificially unstable operating point. A full
nonlinear stability theorem would additionally require well-posedness,
uniqueness, and neighborhood invariance for the closed-loop
hyperbolic--parabolic PDE; robustness and output-feedback extensions
remain natural directions for future work.

\appendix

\bibliographystyle{elsarticle-harv}
\bibliography{biblio}

\end{document}